\documentclass[a4paper,12pt]{article}
 \usepackage{latexsym} 
 \usepackage{amssymb,amsmath} 
\usepackage{cite}
\newtheorem{theorem}{Theorem}
\newtheorem{remark}{Remark}
 \newtheorem{lemma}{Lemma}
 \newtheorem{proposition}{Proposition}
 
 \newtheorem{definition}{Definition} 
\newcommand{\R}{{\mathbb R}}

\begin{document}

\title{Solutions of inhomogeneous perturbed generalized Moisil-Teodorescu system and Maxwell's equations in Euclidean Space}

\author{Juan Bory-Reyes$^{(1)}$ and Marco Antonio P\'{e}rez-de la Rosa$^{(2)}$}

\date{\small $^{(1)}$ ESIME-Zacatenco. Instituto Polit\'ecnico Nacional. CD-MX. 07738. M\'exico.\\E-mail: juanboryreyes@yahoo.com\\
	$^{(2)}$ Department of Actuarial Sciences, Physics and Mathematics, Universidad de las Am\'{e}ricas Puebla.
	San Andr\'{e}s Cholula, Puebla. 72810. M\'{e}xico.\\ Email: marco.perez@udlap.mx}
\maketitle

\begin{abstract} 
\noindent
In this paper, based on a proposed notion of generalized conjugate harmonic pairs in the framework of complex Clifford analysis, necessary and sufficient conditions for the solvability of inhomogeneous perturbed generalized Moisil-Teodorescu systems in higher dimensional Euclidean spaces are proved. As an application, we derive corresponding solvability conditions for the inhomogeneous Maxwell's equations.
\end{abstract}


\section{Introduction}

Clifford analysis is the study of properties of solutions of the first-order, vector-valued Dirac operator $\partial_x$ acting on functions defined on Euclidean spaces $\mathbb R^{m+1}$ ($m\geq 2$) with values in the corresponding real or complex Clifford algebra, that will be denoted below by $\mathbb{R}_{0,m+1}$ and $\mathbb{C}_{0,m+1}$ respectively. Thereby, this function theory may be considered as an elegant way of extending the theory of holomorphic functions in the complex plane to higher dimension and it provides at the same time a refinement of the theory of harmonic functions.

Clifford analysis is centered around the notion of monogenic function, i.e. a null solution of $\partial_x$. It is, however, often important (and interesting) to consider special types of solutions obtained by considering functions taking values in suitable subspaces of the real or complex Clifford algebras (so be the case).

As is established in \cite{BDS2}, there exists a isomorphism between the Cartan algebra of differential forms and the algebra of multivector functions in Clifford analysis. In particular, the action of the operator $d-d^*$, where $d$ and $d^*$ are the differential and codifferential operators (the standard de Rham differential and its adjoint) respectively, on the space of smooth $k$-forms is identified with the action (on the right) of the Dirac operator, which plays the role of the Cauchy-Riemann operator on the space of smooth $k$-vector fields. Meanwhile the action of the operator $d+d^*$ is identified with the action (on the left) of the Dirac operator. A smooth differential form belonging to the kernel of $d+d^*$ was called in \cite {C, Ma} self-conjugate differential form. 

The Moisil-Teodorescu elliptic system of equations of first order in $\mathbb R^{3}$ is a vector valued analogue of Cauchy-Riemann system \cite{MT}.

In the context of real Clifford analysis a generalized Moisil-Teodorescu systems of type $(r,p,q)$ was introduced in \cite {RBDS}, where some general properties of solutions to this system have been investigated. Afterwards, there was a growing interest in the study and better understanding of properties of solutions of generalized Moisil-Teodorescu systems, see for instance \cite {BD1, D1, D2, D3, DLS, FDSc, La, SiHe, So1, So2}. Some special cases of these systems are well known and well understood.

Following the identification mentioned before, a subsystem of generalized Moisil-Teodorescu systems leads to a subsystem of self-conjugate differential forms and vice versa.
 
To deal with the inhomogeneous generalized Moisil-Teodorescu systems in \cite{BoPe} the authors embedded the systems in an appropriate real Clifford analysis setting. Necessary and sufficient conditions for the solvability of inhomogeneous systems are provided and its general solution described and consequently some results in the literature are re-obtained, such as those given in \cite{RBMS, RB, BAPS, PSV1, PSV2}.

The present paper is devoted to give explicit general solution of the inhomogeneous generalized perturbed Moisil-Teodorescu system in the framework of complex Clifford analysis upon the usage of a notion of generalized harmonic conjugates pairs.

A wealth of information about the subject of the nontrivial connection between Maxwell's electrodynamics, Clifford algebras and, in fact, expressing the monogenicity of a certain Clifford algebra-valued function can be found in the literature, see for instance \cite{Chi, F, I, J, K, Ma, McMi, Mi1, Mi2, Mi3, MoWiA, Se, Sp}. The common feature of the method is to represent the Maxwell's equations in a Dirac like form. 

To illustrate the application of the main result, we establish solvability conditions for the inhomogeneous Maxwell's equations. Our study is based on the complex Clifford algebra-based form of the Dirac equation and the Maxwell's equations proposed in \cite{McMi}.

\section{Rudiments of Clifford analysis} 

The section provides a brief exposition of the basic notions and terminology of Clifford analysis aimed at readers who are unfamiliar with this function theory. Standard references are the monographs \cite{BDS1, DSS, GM, GS, GHS}.

Let $\mathbb R^{0,m+1}$ be the real vector space $\mathbb R^{m+1}$ equipped with a quadratic form of signature $(0, m+ 1)$ and
let $e_0, e_1,e_2,\dots, e_m$ be an orthogonal basis of $\mathbb R^{0,m+1}$.

The real Clifford algebra $\R_{0,m+1}$ with generators $e_0, e_1,e_2,\dots, e_m$, subject to the basic multiplication rules 
\[
e_i^2=-1,\quad e_ie_{j}=-e_{j}e_i,\quad i,j=1,2,\dots m,\quad i<j,
\] 
is a real linear associative but non-commutative algebra with identity 1, having dimension $2^{m+1}$ and containing $\R$ and $\R^{m+1}$ as subspaces.

The complex Clifford algebra $\mathbb{C}_{0,m+1}$  constructed over $\mathbb R^{0,m+1}$, as a linear associative algebra over $\mathbb{C}$, has dimension $2^{m+1}$ meaning that one takes the same standard basis as for $\R_{0,m+1}$, with the same multiplication rules, however allowing for complex constants. Indeed, an element of $\mathbb{C}_{0,m+1}$ may be written as $a=\sum_{A} a_A e_A$, where $a_A$ are complex constants and $A$ runs over all the possible ordered sets
\[
A=\{i_{1},\dots,i_{s}\},\quad 0\le i_{1}<i_{2}<\dots<i_{s}\le m,\;{\mbox{or}}\; A=\emptyset,
\]
and
\[
(e_{A}: |A|=s,\; s=0,1,\dots,m+1),\;e_A=e_{i_1}e_{i_2}\cdots e_{i_s},\;e_0=e_\emptyset=1.
\]

One of the basic properties relied upon in building up the $\mathbb{C}_{0, m+1}$-valued continuously differentiable function theory in domains of $\mathbb R^{m+1}$ is the fact that the Dirac operator $\partial_x$ in $\mathbb R^{m+1}$ factorizes the Laplacian $\Delta_x$ through the relation $\partial_x^2=-\Delta_x$, where $\partial_x=\sum_{i=0}^{m}e_i\partial_{x_i}$, $x=(x_0,x_1,\dots ,x_m)\in\mathbb R^{m+1}$.

For  technical reasons to become clear below (see \cite{McMi}), we embed everything into a larger Clifford algebra, say $\mathbb{R}^{0,m+1}\subseteq\mathbb{C}_{0,m+1}\subseteq\mathbb{C}_{0,m+2}$.
Fix $\alpha\in\mathbb{C}$ and set
\[\partial_{x,\alpha}=\partial_x+\alpha e_{m+1},\]
then $-\partial_{x,\alpha}^2=\Delta_x+\alpha^2$, the Helmholtz operator.

If $F$ is a $\mathbb{C}_{0,m+2}$-valued function defined in an open subset $\Omega\subset\mathbb R^{m+1}$, set
\[\partial_{x,\alpha}F:=\sum_{i=0}^{m}e_i\frac{\partial F}{\partial x_i}+\alpha e_{m+1}F,\]
and
\[F \partial_{x,\alpha}:=\sum_{i=0}^{m}\frac{\partial F}{\partial x_i}e_i+\alpha F e_{m+1}.\]

Let $F:\Omega\to\mathbb{C}_{0, m+2}$, whose components are of class $C^1$ in $\Omega$. Then $F$ is called left (right, two-sided) $\alpha$-monogenic in $\Omega$ if $\partial_{x,\alpha}F$ ($F \partial_{x,\alpha}$, or both $\partial_{x,\alpha}F$ and $F \partial_{x,\alpha}$) $=0$ in $\Omega$. Observe that each component of a $\alpha$-monogenic function is annihilated by the Helmholtz
operator $\Delta_x+\alpha^2$. 

An important example of a function which is both right and left $\alpha$-monogenic is the fundamental solution of the perturbed Dirac operator (see \cite{Mi1}), given for $x\in\R^{m+1}\setminus\{0\}$ by
\[E_{\alpha}(x)=
\begin{cases}\displaystyle \frac{1}{\sigma_{m+1}}\frac{\bar{x}}{|x|^{m+1}}-\alpha e_{m+1}\Lambda_{\alpha}(x)+O\left(|x|^{-m+2}\right)&\text{as}\;|x|\to0,\\\\
\displaystyle O\left(\text{exp}\{-\text{Im}\alpha|x|\}\right)&\text{as}\;|x|\to+\infty,
\end{cases}\] 
where $\sigma_{m+1}$ is the area of the sphere in $\R^{m+1}$, $\bar{x}$ is the conjugate of $x$ defined below and for $x\neq0$
\[\Lambda(x):=-\frac{1}{(4\pi)^{\frac{m+1}{2}}}\int_0^{+\infty}\text{exp}\left(\alpha^2 t-\frac{|x|^2}{4t}\right)\frac{dt}{t^{\frac{m+1}{2}}},\]
is the fundamental solution of the Helmholtz operator $\Delta_x+\alpha^2$. The function $E_{\alpha}(x)$ plays the same role in Clifford analysis as the Cauchy kernel does in complex analysis, for this reason it is also called the Cauchy kernel in $\R^{m+1}$.

Writing $E_{\alpha}(x)=E_1(x)+e_{m+1}E_2(x)$ with $E_1(x):=\displaystyle\frac{1}{\sigma_{m+1}}\frac{\bar{x}}{|x|^{m+1}}$ and $E_2(x):=-\alpha \Lambda_{\alpha}(x)$ by letting $|x|\to 0$, we have that $E_1$ is $\mathbb{C}_{0,m+2}^{(1)}$--valued while $E_2$ is $\mathbb{C}_{0,m+2}^{(0)}$--valued. 

If $S$ is a subspace of $\mathbb{C}_{0,m+2}$, then ${\cal E}(\Omega,S)$; ${\cal M}(\Omega,S)$ and ${\cal H}(\Omega,S)$ denote, respectively, the spaces of smooth $S$--valued functions, left monogenic and harmonic $S$--valued functions in $\Omega$. Clearly, we have that ${\cal M}(\Omega,S)\subset{\cal H}(\Omega,S)\subset{\cal E}(\Omega,S)$.

Also, recall that the space $\mathbb{C}_{0,m+2}^{(s)}$ of $s$-vectors in $\mathbb{C}_{0,m+2}$ ($0\le s\le m+2$) is defined by
\begin{equation}\label{1.4}
\mathbb{C}_{0,m+2}^{(s)}={\mbox{span}}_\mathbb{C}(e_A:|A|=s).
\end{equation}
Notice, in particular, that for $s=0$, $\mathbb{C}^{(0)}_{0,m+2}\cong\mathbb{C}$.

For $0\le s\le m+2$ fixed, the space $\mathbb{C}_{0,m+2}^{(s)}$ of $s$-vectors lead to the decomposition
\begin{equation}\label{2.1}
\mathbb{C}_{0,m+2} = \sum^{m+2}_{s=0} \bigoplus \mathbb{C}^{(s)}_{0,m+2},
\end{equation}
and the associated projection operators $[\,\,]_s:\mathbb{C}_{0,m+2}\mapsto\mathbb{C}^{(s)}_{0,m+2}$.

An element $x=(x_0,x_1,\dots,x_{m+1})\in\R^{m+2}$ is usually identified with $x=\sum_{i=0}^{m+1}e_ix_i\in\R^{0,m+2}$.

For $x,y\in\mathbb{C}^{(1)}_{0,m+2}$, the product $xy$ splits in two parts, namely
\begin{equation}\label{2.2}
xy=x\bullet y+x\wedge y,
\end{equation}
where $x\bullet y=[xy]_0$ is the scalar part of $xy$ and $x\wedge y=[xy]_2$ is the 2-vector or bivector part of $xy$ which are given by
\[x\bullet y=-\sum_{i=0}^{m+1}x_iy_i,\]
and
\[x\wedge y=\sum_{i<j}e_ie_j(x_iy_j-x_jy_i).
\]

More generally, for $x\in\mathbb{C}^{(1)}_{0,m+2}$ and $\upsilon\in\mathbb{C}^{(s)}_{0,m+2}$, ($0<s<m+2$), we have that the product $x\upsilon$ decomposes into
\[x\upsilon=x\bullet\upsilon+x\wedge\upsilon,\]
where
\begin{equation}\label{bullet}
x\bullet\upsilon=[x\upsilon]_{s-1}=\frac{1}{2}(x\upsilon-(-1)^s\upsilon x),
\end{equation}
and
\begin{equation}\label{wedge}
x\wedge\upsilon=[x\upsilon]_{s+1}=\frac{1}{2}(x\upsilon+(-1)^s\upsilon x).
\end{equation}

Finally, a conjugation is defined as the unique linear morphism of $\mathbb{C}^{(0)}_{0,m+2}$ with $\bar{e}_0=e_0$, $\bar{e}_j=-e_j$, $j=1,...,m+1$, while for $x,y\in\mathbb{C}^{(0)}_{0,m+2}$,
\[\overline{(xy)}=\bar{y}\,\bar{x}.\]
Notice that for any basic element $e_{A}$ with $|A|=s$, $\bar{e}_{A}=(-1)^{\frac{s(s+1)}{2}}e_{A}.$

\section{Generalized perturbed Moisil-Teodorescu systems}

Let $r, p, q\in\mathbb N$ with $0\le r\le m+2$, $0\le p\le q$ and $r+2q\le m+2$. Exploring further the multivector structure of $\mathbb{C}_{0, m+2}^{(r,p,q)}$ one may also write 
\[\mathbb{C}_{0, m+2}^{(r,p,q)}=\sum_{j=p}^q\bigoplus\mathbb{C}_{0, m+2}^{(r+2j)}.\]

If a $\mathbb{C}_{0, m+2}^{(r+1,p,q)}\bigoplus e_{m+1}\mathbb{C}_{0, m+2}^{(r,p,q)}$--valued smooth function $F$ defined in an open subset $\Omega\subset\mathbb R^{m+1}$ is decomposed following
\begin{align*}
F&=\sum_{j=p}^q\bigoplus F^{(r+2j+1)}+e_{m+1}\sum_{j=p}^q\bigoplus F^{(r+2j)}\\
&=\sum_{j=p}^q\bigoplus F^{(r+2j+1)}+\sum_{j=p}^q\bigoplus (-1)^{r+2j}F^{(r+2j)}e_{m+1},
\end{align*}
then in $\Omega$:
\[\partial_{x,\alpha}F=0\quad\text{if and only if},\]
\begin{equation}\label{MT}
\begin{cases}
\displaystyle \partial_{x}^{-}F^{r+2p+1}-\alpha F^{r+2p}=0,\\
\displaystyle \partial_{x}^{+}F^{r+2j+1}+\partial_{x}^{-}F^{r+2(j+1)+1}-\alpha F^{r+2(j+1)}=0,\qquad j=p,\dots,{q-1};\\
\displaystyle \partial_{x}^{+}F^{r+2q+1}=0,\\
\displaystyle \partial_{x}^{-}(-1)^{r+2p}F^{r+2p}=0,\\
\displaystyle \partial_{x}^{+}(-1)^{r+2j}F^{r+2j}+\partial_{x}^{-}(-1)^{r+2(j+1)}F^{r+2(j+1)}+\alpha (-1)^{r+2j+1}F^{r+2j+1}=0,\\
\displaystyle \hfill j=p,\dots,{q-1};\\
\displaystyle \partial_{x}^{+}(-1)^{2+2q}F^{r+2q}+\alpha (-1)^{r+2q+1}F^{r+2q+1}=0,
\end{cases}
\end{equation}
where the differential operators $\partial_x^+$ and $\partial_x^-$ act on smooth $\mathbb{C}_{0,m+2}^{(s)}$--valued functions $F^s$ in $\Omega$ as
\begin{equation}\label{d+}
\partial_x^+F^s=\frac{1}{2}(\partial_xF^s-(-1)^sF^s\partial_x),
\end{equation}
and
\begin{equation}\label{d-}
\partial_x^-F^s=\frac{1}{2}(\partial_xF^s+(-1)^sF^s\partial_x).
\end{equation}
It is perhaps worth remarking that $\partial^+_{x} F^s$ is $\mathbb{C}^{(s+1)}_{0,m+2}$--valued while $\partial^-_{x} F^s$ is $\mathbb{C}^{(s-1)}_{0,m+2}$--valued.

The system (\ref{MT}) generalizes that of \cite{BoPe} and is called generalized perturbed Moisil-Teodorescu system of type $(r,p,q)$. 

In \cite{BoPe} is pointed out that for $\alpha=0$ the system (\ref{MT}) includes some basic systems of first order linear partial differential equations as particular cases. For example, if $p=0$, $q=1$ and $0\le r\le m+1$ fixed, the system (\ref{MT}) reduces to the Moisil-Teodorescu system in $\R^{m+1}$ introduced  in \cite{BD1}. If $p=q=0$ and $0<r<m+1$ fixed, the system (\ref{MT}) reduces to the generalized Riesz system $\partial_xF^r=0$; its solutions are called harmonic multi-vector fields. If $r=0$, $p=0$ and $m+1=3$ then $q=1$, the original Moisil-Teodorescu system introduced in \cite{Shap} is re-obtained. If $r=0$, $p=0$ and $m+1=4$ then $q=2$ and one obtains the Fueter system in $\R^4$ for so-called left regular functions of quaternion variable; it lies at the basis of quaternionic analysis (see \cite{Fu, Sub}).

\section{The inhomogeneous Dirac equation}

From now on we assume $\Omega$ to be a Lipschitz domain in $\R^{m+1}$, i.e., a domain whose boundary $\Gamma$ is given locally by the graph of a real valued Lipschitz function, after an appropriate rotation of coordinates.

The fundamental tool for solving the inhomogeneous perturbed Dirac equation (commonly called the $\overline\partial$-problem)
\begin{equation}\label{inh sys Cliff}
	\partial_{x,\alpha} F=G,
	\end{equation}
is the Borel-Pompeiu integral formula (see below) which is named after the French and Romanian mathematicians \'Emile Borel (1871-1956) and Dimitrie Pompeiu (1873-1954), respectively.

For bounded $F\in C^0(\Omega;\mathbb{C}_{0,m+2})$, we consider the Teodorescu and the Cauchy-type operators associated to the Cauchy
kernel $E_{\alpha}$, i.e., 
\[T_{\Omega}[F](x):=\int_{\Omega}E_{\alpha}(x-y)F(y)dy,\quad x\in\R^{m+2},\]
and by
\[C_{\Gamma}[F](x):=-\int_{\Gamma}E_{\alpha}(x-y)n(y)F(y)d\Gamma_y,\quad x\notin\Gamma,\]
where $n(y)=\sum_{i=0}^{m+1}e_i n_i(y)$ is the outward pointing unit normal to $\Gamma$ at $y\in\Gamma$.

\begin{lemma}
	Let $F\in C^1(\Omega;\mathbb{C}_{0,m+2})\cap C^0(\Omega\cup \Gamma;\mathbb{C}_{0,m+2})$. Then we have
	\begin{itemize}
		\item [i)] \begin{equation}\label{BPCliff}
		C_{\Gamma}[F](x)+T_{\Omega}[\partial_{x,\alpha} F](x)=\begin{cases}
		F(x),&x\in\Omega,\\
		0,&x\in\R^{m+1}\setminus(\Omega\cup \Gamma).
		\end{cases}
		\end{equation}
		\item [ii)] \begin{equation}\label{RightInv}
		\partial_{x,\alpha} T_{\Omega}[F](x)=\begin{cases}
		F(x),&x\in\Omega,\\
		0,&x\in\R^{m+1}\setminus(\Omega\cup \Gamma).
		\end{cases}
		\end{equation}
		\item [iii)] \begin{equation}
		\partial_{x,\alpha} C_{\Gamma}[F](x)=0,\quad x\in\Omega\cup\left(\R^{m+1}\setminus(\Omega\cup \Gamma)\right).
		\end{equation}
	\end{itemize}
\end{lemma}

For a proof of the Borel-Pompeiu formula (\ref{BPCliff}) see, e.g. \cite{Zhen}.

\begin{remark}
The Borel-Pompeiu formula (\ref{BPCliff}) solves the inhomogeneous perturbed Dirac equation (\ref{inh sys Cliff}) in the standard way and the general solution is given by
	\begin{equation}\label{inh sys sol Cliff}
	F=T_{\Omega}[G]+H,
	\end{equation}
where $H\in\mathcal M(\Omega; \mathbb{C}_{0,m+2})$.
\end{remark}

\section{Generalized conjugate harmonics pairs}

The notion of conjugate harmonic functions in the complex plane is well-known. This concept has been generalized to higher dimensional setting in the framework of Clifford analysis, see for instance \cite{BoPe,FDSo, GM, Nol, Shap}.

In this section, we introduce a new generalization of notion of conjugate harmonic functions in a Clifford setting, based on a certain splitting of the Clifford algebra.
\begin{definition}
Let \[F_1=\sum_{j=p}^q F^{r+2j+1}+e_{m+1}\sum_{j=p}^q F^{r+2j},\]
in ${\cal H}\left(\Omega;\mathbb{C}_{0, m+2}^{(r+1,p,q)}\bigoplus e_{m+1}\mathbb{C}_{0, m+2}^{(r,p,q)}\right)$. An \[F_2=\left(F^{r+2p-1}+F^{r+2q+3}\right)+e_{m+1}\left(F^{r+2p-2}+F^{r+2q+2}\right),\] 
in ${\cal H}\left(\Omega;\left(\mathbb{C}^{(r+2p-1)}_{0,m+2} \bigoplus \mathbb{C}^{(r+2q+3)}_{0,m+2}\right)\bigoplus e_{m+1}\left(\mathbb{C}^{(r+2p-2)}_{0,m+2} \bigoplus \mathbb{C}^{(r+2q+2)}_{0,m+2}\right)\right)$ is called hyper-conjugate harmonic to $F_1$ if 
\[F_1+F_2\in {\cal M}(\Omega;\mathbb{C}_{0,m+2}).\]
The pair $(F_1,F_2)$ is then called a pair of hyper-conjugate harmonic functions.
\end{definition}

\section{Main result}

We are in condition to state and proof our main result
\begin{theorem}\label{maintheo}
Let $\displaystyle G\in C\left(\Omega;\sum_{j=p}^{q+1} \bigoplus \mathbb{C}^{(r+2j)}_{0,m+2}+\sum_{j=p}^{q+1} \bigoplus \mathbb{C}^{(r+2j-1)}_{0,m+2}e_{m+1}\right)$. The inhomogeneous generalized perturbed Moisil-Teodorescu system

\begin{equation}\label{general system}
\begin{cases}
 \displaystyle \partial_{x}^{-}F^{r+2p+1}-\alpha F^{r+2p}=G^{r+2p},\\
 \displaystyle \partial_{x}^{+}F^{r+2j+1}+\partial_{x}^{-}F^{r+2(j+1)+1}-\alpha F^{r+2(j+1)}=G^{r+2j+2},\qquad j=p,\dots,{q-1};\\
 \displaystyle \partial_{x}^{+}F^{r+2q+1}=G^{r+2q+2},\\
 \displaystyle \partial_{x}^{-}(-1)^{r+2p}F^{r+2p}=G^{r+2p-1},\\
 \displaystyle \partial_{x}^{+}(-1)^{r+2j}F^{r+2j}+\partial_{x}^{-}(-1)^{r+2(j+1)}F^{r+2(j+1)}+\\
 \displaystyle\quad+\alpha (-1)^{r+2j+1}F^{r+2j+1}=G^{r+2j+1}, \hfill j=p,\dots,{q-1};\\
 \displaystyle \partial_{x}^{+}(-1)^{2+2q}F^{r+2q}+\alpha (-1)^{r+2q+1}F^{r+2q+1}=G^{r+2q+1},
\end{cases}
\end{equation}
where $G^s\in C^1\left(\Omega;\mathbb{C}^{(s)}_{0,m+2}\right)\cap C^0\left(\Omega\cup \Gamma;\mathbb{C}^{(s)}_{0,m+2}\right)$, has a solution if and only if for the $\left(\mathbb{C}^{(r+2p-1)}_{0,m+2} \bigoplus\mathbb{C}^{(r+2q+3)}_{0,m+2}\right)\bigoplus e_{m+1}\left(\mathbb{C}^{(r+2p-2)}_{0,m+2} \bigoplus \mathbb{C}^{(r+2q+2)}_{0,m+2}\right)$--valued function
	\begin{align*}
	 P:&=\left(\int E_1\bullet G^{r+2p}dy+(-1)^{r+2p}\int E_2\,G^{r+2p-1}dy+\right.\\
	 &\quad\left.+\int E_1\wedge G^{r+2q+2}dy\right)+\\
	&\quad+e_{m+1}\left((-1)^{r+2p}\int E_1\bullet G^{r+2p-1}dy+\int E_1\wedge G^{r+2q+1}dy+\right.\\
	&\quad\left.+\int E_2\,G^{r+2q+2}dy\right),
	\end{align*}
	\begin{itemize}
		\item[(A)] either is identically zero;
		\item[(B)] or has a hyper-conjugate harmonic function.
	\end{itemize}
	If it is true, then the general solution of (\ref{general system}) is given by:
	\begin{itemize}
		\item[(A*)] either
	\end{itemize}
		\begin{align}
	F=&\sum_{j=p}^{q}\left[\left(\int E_1\wedge G^{r+2j}dy+\int E_1\bullet G^{r+2j+2}dy+\right.\right.\notag\\
	&+\left.(-1)^{r+2j+2}\int E_2\,G^{r+2j+1}dy\right)+\notag\\
	&+e_{m+1}\left(\int E_1\wedge G^{r+2j-1}dy+(-1)^{r+2j+2}\int E_1\bullet G^{r+2j+1}dy+\right.\notag\\
	&\left.\left.+\int E_2\,G^{r+2j}dy\right)\right] +\hat{H}_1,\label{sol A*}
		\end{align}
	\begin{itemize}
		\item[(B*)] or
	\end{itemize}
		\begin{align}
		 F=&\sum_{j=p}^{q}\left[\left(\int E_1\wedge G^{r+2j}dy+\int E_1\bullet G^{r+2j+2}dy+\right.\right.\notag\\
		&+\left.(-1)^{r+2j+2}\int E_2\,G^{r+2j+1}dy\right)+\notag\\
		&+e_{m+1}\left(\int E_1\wedge G^{r+2j-1}dy+(-1)^{r+2j+2}\int E_1\bullet G^{r+2j+1}dy+\right.\notag\\
		&\left.\left.+\int E_2\,G^{r+2j}dy\right)\right] +\hat{H}_1+\tilde{H}_1,\label{sol B*}
		\end{align}
where $\tilde{H}_1$ is a harmonic hyper-conjugate of $-P$, and $\hat{H}_1$ is an arbitrary monogenic function.	
\end{theorem}
{\bf Proof:}
First, notice that system (\ref{general system}) is a restriction of (\ref{inh sys Cliff}), for 
\[F=\sum_{j=p}^q F^{r+2j+1}+e_{m+1}\sum_{j=p}^q F^{r+2j},\] 
and with 
\[G=\sum_{j=p}^{q+1} G^{r+2j}+\sum_{j=p}^{q+1} G^{r+2j-1}e_{m+1}.\] 
Let $F$ be such a solution. Looking at the components of (\ref{inh sys sol Cliff}) one has:
\begin{equation}\label{system sol ICS}
\begin{cases}
\displaystyle 0=\int E_1\bullet G^{r+2p}dy+(-1)^{r+2p}\int E_2\,G^{r+2p-1}dy+H^{r+2p-1},\\
\displaystyle F^{r+2j+1}=\int E_1\wedge G^{r+2j}dy+\int E_1\bullet G^{r+2j+2}dy+\\
\displaystyle\qquad\qquad+(-1)^{r+2j+2}\int E_2\,G^{r+2j+1}dy+H^{r+2j+1},\qquad j=p,...,q,\\
\displaystyle 0=\int E_1\wedge G^{r+2q+2}dy+H^{r+2q+3},\\

\displaystyle 0=(-1)^{r+2p}\int E_1\bullet G^{r+2p-1}dy+H^{r+2p-2},\\
\displaystyle F^{r+2j}=\int E_1\wedge G^{r+2j-1}dy+(-1)^{r+2j+2}\int E_1\bullet G^{r+2j+1}dy+\\
\displaystyle\qquad\qquad+\int E_2\,G^{r+2j}dy+H^{r+2j},\hfill j=p,...,q,\\
\displaystyle 0=\int E_1\wedge G^{r+2q+1}dy+\int E_2\,G^{r+2q+2}+H^{r+2q+2},
\end{cases}
\end{equation}
where $(H_1,H_2)$ is a hyper-conjugate harmonic pair with \[H_1:=\sum_{j=p}^q H^{r+2j+1}+e_{m+1}\sum_{j=p}^q H^{r+2j},\]
and 
\[H_2:=\left(H^{r+2p-1}+H^{r+2q+3}\right)+e_{m+1}\left(H^{r+2p-2}+H^{r+2q+2}\right).\] 
For the case of $P$ being the zero function, one obtains that $\hat{H}_1$ becomes also the zero function, and both the necessity of (A) and the formula (\ref{sol A*}) is proved. 

For the case of $P$ being identically zero one has
\begin{align*}
H_2:&=-\left(\int E_1\bullet G^{r+2p}dy+(-1)^{r+2p}\int E_2\,G^{r+2p-1}dy+\right.\\
&\quad\left.+\int E_1\wedge G^{r+2q+2}dy\right)-\\
&\quad-e_{m+1}\left((-1)^{r+2p}\int E_1\bullet G^{r+2p-1}dy+\int E_1\wedge G^{r+2q+1}dy+\right.\\
&\quad\left.+\int E_2\,G^{r+2q+2}dy\right),
\end{align*} 
meaning the existence of the hyper-conjugate harmonic function to 
\begin{align*}
&-\left(\int E_1\bullet G^{r+2p}dy+(-1)^{r+2p}\int E_2\,G^{r+2p-1}dy+\right.\\
&+\left.\int E_1\wedge G^{r+2q+2}dy\right)-\\
&-e_{m+1}\left((-1)^{r+2p}\int E_1\bullet G^{r+2p-1}dy+\int E_1\wedge G^{r+2q+1}dy+\right.\\
&\left.+\int E_2\,G^{r+2q+2}dy\right),
\end{align*} 
which (whenever it exists) will be denote by $\tilde{H}_1$. Thus, the necessity of (B) is also proved together with formula (\ref{sol B*}).

Finally, for the ``sufficiency part'' one has to reverse the reasoning. If (A) is true the one can take $H_2$ as the zero function implying $H_1$ to be an arbitrary monogenic function; thus we arrive at (\ref{sol A*}). If (B) is true then one can take in (\ref{system sol ICS}) the function $H_2$ to be $H_2=-P$; hence, one can write $H_1$ as $H_1=\hat{H}_1+\tilde{H}_1$ for any arbitrary monogenic $\hat{H}_1$; thus, one obtains (\ref{sol B*}), and the theorem is proved.
$\hfill \square$

\begin{remark}
Our approach provides a rather new generalization and certain consolidation of the hyper-complex method (by means of techniques based on harmonic conjugates) described in \cite{RBMS, RB, BAPS, BoPe, CLSSS, DP, DK, PSV1, PSV2} concerning the solvability of the $\overline\partial$-problem in spaces of functions satisfying different systems in the framework of (complex) Clifford and Quaternionic analysis. 
\end{remark}

\section{The inhomogenous Maxwell's equations}
To keep the exposition self-contained let us repeat key observations from \cite[pag. 1613]{McMi}. Let $E$ be a $r$--vector and $H$ a $(r+1)$--vector with $0\leq r\leq m+1$, defined in some open subset of $\mathbb{R}^{m+1}$. To these, we associate an $\mathbb{C}_{0,m+2}$--valued function $M$ by setting
\begin{equation}\label{EHfield}
M:=H-ie_{m+1}E=H+i(-1)^{r+1}Ee_{m+1}.
\end{equation}

\begin{remark}
Observe that
\[\partial_{x,\alpha}M=i(-1)^{r+1}\left(\partial_{x}^{+}E+\partial_{x}^{-}E-i\alpha H\right)e_{m+1}+\left(\partial_{x}^{+}H+\partial_{x}^{-}H+i\alpha E\right).\]
\end{remark}

\begin{proposition}
The function $M$ defined in (\ref{EHfield}) is (left or right) $\alpha$-monogenic if and only if $E$ and $H$ satisfy the Maxwell's equations
\begin{equation}\label{Max}
\begin{cases}
\displaystyle \partial_{x}^{+}E-i\alpha H=0,\\
\displaystyle \partial_{x}^{-}E=0,\\
\displaystyle \partial_{x}^{-}H+i\alpha E=0,\\
\displaystyle \partial_{x}^{+}H=0.
\end{cases}
\end{equation}
Moreover, for $x\in\mathbb{R}^{m+1}\setminus\{0\}$ the function $M$ satisfies the radiation condition
\[\left(1-ie_{m+1}\frac{x}{|x|}\right)M(x)=o\left(|x|^{-m/2}\right)\quad\text{as}\quad |x|\to\infty,\]
if and only if $E$ and $H$ satisfy the Silver-M\"{u}ller-type radiation conditions
\begin{align}
\frac{x}{|x|}\wedge E-H&=o\left(|x|^{-m/2}\right)\quad\text{as}\quad |x|\to\infty,\\
\frac{x}{|x|}\bullet H-E&=o\left(|x|^{-m/2}\right)\quad\text{as}\quad |x|\to\infty.
\end{align}
\end{proposition}

\begin{remark}
If $\alpha$ is non-zero then, because $(\partial_{x}^{+})^2=(\partial_{x}^{-})^2=0$,  the equations $\partial_{x}^{-}E=0$ and $\partial_{x}^{+}H=0$ become superfluous. Nonetheless, as for $\alpha=0$ Maxwell's equations decouple (i.e. $E$ and $H$ become unrelated), it is precisely this case for which these two equations are relevant. Also, when $m+1=3$, $r=1$ (and $\alpha\neq0$), the formulae (\ref{Max}) reduce to the
more familiar system of equations
\begin{equation}\label{Max2}
\begin{cases}
	\displaystyle \nabla\times E-i\alpha H=0,\\
	\displaystyle \nabla\times H+i\alpha E=0.
\end{cases}
\end{equation}
\end{remark}

\begin{theorem}\label{maintheoMax}
	Let $\displaystyle G\in C\left(\Omega;\sum_{j=0}^{1} \bigoplus \mathbb{C}^{(r+2j)}_{0,m+2}+\sum_{j=0}^{1} \bigoplus \mathbb{C}^{(r+2j-1)}_{0,m+2}e_{m+1}\right)$. The inhomogeneous Maxwell's equations
	
	\begin{equation}\label{InhMax}
	\begin{cases}
	\displaystyle \partial_{x}^{+}E-i\alpha H=G^{r+1},\\
	\displaystyle \partial_{x}^{-}E=G^{r-1},\\
	\displaystyle \partial_{x}^{-}H+i\alpha E=G^{r},\\
	\displaystyle \partial_{x}^{+}H=G^{r+2},
	\end{cases}
	\end{equation}
	where $G^s\in C^1\left(\Omega;\mathbb{C}^{(s)}_{0,m+2}\right)\cap C^0\left(\Omega\cup \Gamma;\mathbb{C}^{(s)}_{0,m+2}\right)$ has a solution if and only if for the $\left(\mathbb{C}^{(r-1)}_{0,m+2} \bigoplus \mathbb{C}^{(r+3)}_{0,m+2}\right)\bigoplus e_{m+1}\left(\mathbb{C}^{(r-2)}_{0,m+2} \bigoplus \mathbb{C}^{(r+2)}_{0,m+2}\right)$--valued function
	\begin{align*}
	P:&=\left(\int E_1\bullet G^{r}dy+(-1)^{r}\int E_2\,G^{r-1}dy+\int E_1\wedge G^{r+2}dy\right)+\\
	&+e_{m+1}\left((-1)^{r}\int E_1\bullet G^{r-1}dy+\int E_1\wedge G^{r+1}dy+\int E_2\,G^{r+2}dy\right),
	\end{align*}
	\begin{itemize}
		\item[(A)] either is identically zero;
		\item[(B)] or has a hyper-conjugate harmonic function.
	\end{itemize}
	If it is true, then the general solution of (\ref{InhMax}) is given by:
	\begin{itemize}
		\item[(A*)] either
		\begin{align}
		M&=H-ie_{m+1}E=H+i(-1)^{r+1}Ee_{m+1}\notag\\
		&=\left(\int E_1\wedge G^{r}dy+\int E_1\bullet G^{r+2}dy+(-1)^{r+2}\int E_2\,G^{r+1}dy\right)+\notag\\
		&+e_{m+1}\left(\int E_1\wedge G^{r-1}dy+(-1)^{r+2}\int E_1\bullet G^{r+1}dy+\int E_2\,G^{r}dy\right)+\notag\\
		&+\hat{H}_1,\label{sol A*2}
		\end{align}
		\item[(B*)] or
		\begin{align}
		M&=H-ie_{m+1}E=H+i(-1)^{r+1}Ee_{m+1}\notag\\
		&=\left(\int E_1\wedge G^{r}dy+\int E_1\bullet G^{r+2}dy+(-1)^{r+2}\int E_2\,G^{r+1}dy\right)+\notag\\
		&+e_{m+1}\left(\int E_1\wedge G^{r-1}dy+(-1)^{r+2}\int E_1\bullet G^{r+1}dy+\int E_2\,G^{r}dy\right)+\notag\\
		&+\hat{H}_1+\tilde{H}_1,\label{sol B*2}
		\end{align}
	\end{itemize}
	where $\tilde{H}_1$ is a harmonic hyper-conjugate of $-P$, and $\hat{H}_1$ is an arbitrary $\alpha$-monogenic function.	
\end{theorem}
{\bf Proof:} It is a direct consequence of the Theorem \ref{maintheo} by choosing $p=q=0$ and $F=M$.
$\hfill \square$

\begin{remark}
Observe that, under the assumptions of Theorem 1, it is possible to obtain both the electric field and the magnetic field separately, say if $P$ is identically zero, then
\begin{align*}
H&=\int E_1\wedge G^{r}dy+\int E_1\bullet G^{r+2}dy+(-1)^{r+2}\int E_2\,G^{r+1}dy+\left[\hat{H}_1\right]_{r+1},\\
E&=i\left(\int E_1\wedge G^{r-1}dy+(-1)^{r+2}\int E_1\bullet G^{r+1}dy+\int E_2\,G^{r}dy+\left[\hat{H}_1\right]_{r}\right),
\end{align*}
while if $P$ has a hyper-conjugate harmonic function, then
\begin{align*}
H&=\int E_1\wedge G^{r}dy+\int E_1\bullet G^{r+2}dy+(-1)^{r+2}\int E_2\,G^{r+1}dy+\\
&\quad+\left[\hat{H}_1\right]_{r+1}+\left[\tilde{H}_1\right]_{r+1},\\
E&=i\left(\int E_1\wedge G^{r-1}dy+(-1)^{r+2}\int E_1\bullet G^{r+1}dy+\int E_2\,G^{r}dy+\right.\\
&\quad+\left.\left[\hat{H}_1\right]_{r}+\left[\tilde{H}_1\right]_{r}\right).
\end{align*}
\end{remark}

\section*{Acknowledgments}
The authors wish to thank Instituto Polit\'ecnico Nacional and Universidad de las Am\'ericas Puebla, for partial financial support.

\end{document}